\documentclass{amsart}
\usepackage[breaklinks,pdfstartview=FitH,colorlinks,linktocpage,bookmarks=true]{hyperref}

\usepackage{graphicx}

\usepackage[cmtip,all]{xy}
\usepackage{upref}

\def\procsection#1{\section{#1}}
\def\procsubsection#1{\subsection*{#1}}
\newenvironment{procfigure}{\begin{figure}}{\end{figure}}
\newenvironment{procthebibliography}[1]{\begin{thebibliography}{#1}}{\end{thebibliography}}

\newtheorem{proc-theorem}{Theorem}[section]
\newtheorem{proc-lemma}[proc-theorem]{Lemma}

\newtheorem{proc-corollary}[proc-theorem]{Corollary}
\newtheorem{proc-claim}[proc-theorem]{Claim}

\theoremstyle{definition}
\newtheorem{proc-definition}[proc-theorem]{Definition}
\newtheorem{proc-example}[proc-theorem]{Example}

\theoremstyle{remark}
\newtheorem{proc-remark}[proc-theorem]{Remark}

\numberwithin{equation}{section}

\DeclareMathOperator{\projection}{pr}
\DeclareMathOperator{\rank}{rk}
\DeclareMathOperator{\relatively}{rel}
\DeclareMathOperator{\interior}{Int}
\DeclareMathOperator{\closure}{Cl}
\DeclareMathOperator{\identity}{Id}
\DeclareMathOperator{\image}{Im}
\DeclareMathOperator{\kernel}{Ker}

\begin{document}

\title{Suspension theorems for links and link maps}
\author{Mikhail Skopenkov}
\address{Department of Differential Geometry\\ Faculty of Mechanics and Mathematics \\
Moscow State University \\ 119992, Moscow, Russia.}
\curraddr{}
\email{skopenkov@rambler.ru, http://skopenkov.ru}
\thanks{The author was supported in part by INTAS grant 06-1000014-6277,
Russian Foundation of Basic Research grants 05-01-00993-a,
06-01-72551-NCNIL-a, 07-01-00648-a, President of the Russian
Federation grant NSh-4578.2006.1, Agency for Education and Science
grant RNP-2.1.1.7988, and Moebius Contest Foundation for Young Scientists.}



\date{}


\commby{Alexander Dranishnikov}

\begin{abstract}
We present a new short proof of the explicit formula for the group of
links (and also link maps) in the 'quadruple point free' dimension.
Denote by $L^m_{p,q}$ (respectively, $C^{m-p}_p$) the group of smooth embeddings
$S^p\sqcup S^q\to S^m$ (respectively, $S^p\to S^m$) up to smooth isotopy.
Denote by $LM^m_{p,q}$ the group of link maps $S^p\sqcup S^q\to S^m$ up to link homotopy.

\smallskip
\noindent\textbf{Theorem~1.} \textit{
If $p\le q\le m-3$ and $2p+2q\le 3m-6$ then}
\begin{equation*}
L^m_{p,q}\cong
\pi_p(S^{m-q-1})\oplus\pi_{p+q+2-m}(SO/SO_{m-p-1})\oplus C^{m-p}_p\oplus
C^{m-q}_q.
\end{equation*}

\noindent\textbf{Theorem~2.} \textit{
If $p, q\le m-3$ and $2p+2q\le 3m-5$ then $LM^m_{p,q}\cong \pi^S_{p+q+1-m}$.}

\smallskip

Our approach is based on the use of the suspension operation for
links and link maps, and \textit{suspension theorems} for them.
\end{abstract}

\maketitle


\procsection{Introduction}\label{proc-sect1}

This paper\footnote{This is an improved version of the paper in Proc. Amer. Math. Soc. 137:1 (2009), 359--369.}
is on knotting problem of higher dimensional manifolds
(for recent surveys see \cite{RS99, Sko07, S}).
We study knots and links in codimension at least 3, where
a complete answer can sometimes be obtained, in contrast to the
classical situation of simple closed curves in $\Bbb R^3$ (see a survey in \cite{M}).

Denote by $L^m_{p,q}$ the set of smooth
embeddings $S^p\sqcup S^q\to S^m$ up to
smooth isotopy. Denote by $C^{m-p}_p$ the set of smooth
embeddings $S^p\to S^m$ up to smooth isotopy.
For $p,q\le m-3$ these sets are commutative groups with respect to 'componentwise connected sum'
operation \cite{Hae66C}, cf.~\cite[Remark~2.3.ab]{Sko06}.

The main result of this paper is a new short proof of an explicit formula for the group
$L_{p,q}^m$ in terms of the groups~ $C^{m-p}_p$, $C^{m-q}_q$ and certain homotopy groups:

\begin{proc-theorem}\label{proc-th1} If $p\le q\le m-3$ and $2p+2q\le 3m-6$ then
\begin{equation*}
L^m_{p,q}\cong\pi_p(S^{m-q-1})\oplus\pi_{p+q+2-m}(V_{M+m-p-1,M})\oplus
C^{m-p}_p\oplus C^{m-q}_q.
\end{equation*}
\end{proc-theorem}

Here $V_{M+l,M}$ is the {\it Stiefel manifold} of $M$-frames at the
origin of $\Bbb R^{M+l}$, where $M$ is large. Many of the groups
$\pi_n(V_{M+l,M})$ and $C^{m-p}_p$ are known \cite{Pae56, Hae66A}.

\begin{proc-example} $L^6_{3,3}\cong \mathbb{Z}\oplus \mathbb{Z}\oplus \mathbb{Z}\oplus \mathbb{Z}$.
\end{proc-example}

Theorem~\ref{proc-th1} is the strongest available {\it readily calculable} classification of $2$-component links in spheres in the sense of \cite[Remark~1.1]{S}. However, for arbitrary $p,q\le m-3$ there are a \emph{rational} classification (see \cite[Theorem~1.9]{CFS14}) and a famous
exact sequence involving the groups $L^m_{p,q}$, certain homotopy
groups and maps between them involving Whitehead products (see \cite[Theorem~1.1]{Hae66C} and \cite{Ha86}).

Theorem~\ref{proc-th1} was proved in \cite[Theorems 10.7 and 2.4]{Hae66C} under
stronger restrictions $p\le q$ and $p+3q\le 3m-7$. However, the Haefliger
argument can be extended to cover the dimension range $2p+2q\le 3m-7$ but not the 'boundary' case $2p+2q= 3m-6$; see Remark~\ref{proc-rem2}. 
The second inequality in~\ref{proc-th1} is sharp; see Remark~\ref{proc-rem1}.


We reduce the classification of links to the classification of link
maps, which is an interesting problem in itself \cite{Sco68, Kos90,
HaKa98}.

{\it A link map} is a continuous map $f:X\sqcup Y\to Z$ such that
$fX\cap fY=\emptyset$. {\it A link homotopy} is a continuous
family of link maps $f_t:X\sqcup Y\to Z$. Denote by $LM^m_{p,q}$ the
set of link maps $S^p\sqcup S^q\to S^m$ up to link homotopy. For
$p,q\le m-3$ this set is a commutative group
with respect to 'componentwise connected sum' operation (by
\cite[p. 187]{Sco68}, \cite[Remark~2.4]{Kos88} and 'link concordance
implies link homotopy' theorem discussed below in this section).

The second result of this paper is a short proof of the following
theorem:

\begin{proc-theorem}\label{proc-th2}\cite{HaKa98} If $p,q\le m-3$ 
and $2p+2q\le 3m-5$ then
$
LM^m_{p,q}\cong \pi^S_{p+q-m+1}.
$
\end{proc-theorem}

The isomorphism is the $\alpha$-invariant (see \S\ref{proc-sect3}).
The second inequality is sharp \cite{HaKa98}. 

Theorem~\ref{proc-th2} is the strongest known readily calculable classification of link
maps for $p,q\le m-3$. However, under slightly weaker
dimension restriction there is an exact sequence involving
the groups $LM^m_{p,q}$ and certain bordism groups \cite[Theorem~A]{Kos90}, cf.~Theorem~\ref{proc-th6} below.

Our approach is based on the use of the suspension map.
{\it The suspension map}
$
\Sigma:LM^m_{p,q}\to LM^{m+1}_{p+1,q}
$ 
is defined by suspending the $p$-component and including the
$q$-component. It is easy to see that for
$M$ large $LM^{m+M}_{p+M,q}\cong\pi^S_{p+q-m+1}$ \cite{Kos88}. Thus Theorem~\ref{proc-th2} follows from the
following assertion:

\begin{proc-theorem}[Suspension theorem for link maps]\label{proc-th2pr}\cite{HaKa98}
If $p,q\le m-3$ then the suspension map is bijective for $2p+2q\le
3m-5$ and surjective for $2p+2q\le 3m-4$.
\end{proc-theorem}

This theorem has been known earlier only as a corollary of
Theorem~\ref{proc-th2}. We give a short direct proof of
Theorem~\ref{proc-th2pr} analogous to Zeeman's proof of the higher-dimensional Poincar\'e
conjecture and using a version of Alexander's trick. Our proof is almost self-contained, we use only 'concordance implies isotopy' theorem and its version for link concordances. In the proof of Theorem~\ref{proc-th1} we use basics of immersion theory and suspension sequences by U.~Koschorke and A.~Skopenkov (cf. the sequence by V.~Nezhinsky \cite{Ne84}).


Let us introduce some notions and conventions.

An embedding $f:X\times I\to S^m\times I$ is {\it a concordance} if
$X\times 0=f^{-1}(S^m\times 0)$ and $X\times 1=f^{-1}(S^m\times 1)$.
We tacitly use the facts that in codimension at least $3$ {\it
concordance implies isotopy} and {\it any concordance or isotopy is
ambient} \cite{Hud69}.

Similarly, {\it a link concordance} is a continuous map $f:(X\sqcup
Y)\times I\to S^m\times I$ such that $f(X\times I)\cap f(Y\times
I)=\emptyset$, $(X\sqcup Y)\times 0=f^{-1}(S^m\times 0)$ and
$(X\sqcup Y)\times 1=f^{-1}(S^m\times 1)$.
We say that a link map $f:S^p\sqcup S^q\to S^m$ is {\it null link concordant},
if it extends to a link map $D^{p+1}\sqcup D^{q+1}\to D^{m+1}$. The latter link map is called {\it null link concordance}.

In codimension at least
$3$ {\it link concordance implies link homotopy,} which was announced in \cite{Mel00}, cf. \cite{BaTe99, KrTa97, Kos97}, and proved in \cite{Mel} (unfortunately, the latter preprint has not been published). This result is essentially used only in the proof of the injectivity in~\ref{proc-th2pr} modulo~\ref{proc-lem3pr}. Without it all the results and proofs of the present paper remain true, but $LM^m_{p,q}$ should be understood as the group of link maps up to link \emph{concordance}.

Other minor variations of known theorems which we use without exact references are Theorem~\ref{proc-embth} and basics of immersion theory in the proof of Theorem~\ref{proc-th6}.


This paper is organized as follows. In \S\ref{proc-sect2} we prove Theorem~\ref{proc-th2pr}. In \S\ref{proc-sect3} we deduce
Theorem~\ref{proc-th1} from Theorem~\ref{proc-th2pr}. Sections \S\ref{proc-sect2} and \S\ref{proc-sect3} can be read
independently from each other.


In \cite{CRS07,CRS} (cf. \cite{Sko06}) a similar approach is applied to the classification of
embeddings $S^p\times S^q\to S^m$.

\procsection{Classification of link maps}\label{proc-sect2}

We prove Theorem~\ref{proc-th2pr} as follows. First we prove the surjectivity in
case $p\le q$. Then we prove analogously the injectivity in case $p\le q$, and finally
we deduce the case $p>q$ of Theorem~\ref{proc-th2pr} from the case $p\le q$.

Let us introduce our main notion and state our main lemma.


\begin{proc-definition}\label{proc-defslm} (see Figure~\ref{proc-fig1})
Let $S^{k}=D^{k}_+\cup (S^{k-1}\times I)\cup D^{k}_-$ be the
standard decomposition of the sphere, where $\partial D^{k}_+=S^{k-1}\times 0=S^{k-1}$ is the equator of~$S^k$.
A link map $f:S^{p}\sqcup S^{q}\to S^{m}$ is {\it standardized} if the
following 3 conditions hold:

\begin{enumerate}
\item $fD^{p}_+\subset D^{m}_+$, \quad $fD^{p}_-\subset D^{m}_-$, \quad
$f(S^{p-1}\times I)\subset S^{m-1}\times I$;\label{proc-item1}

\item $fS^q\subset S^{m-1}\times I$;\label{proc-item2}

\item $f(S^{p-1}\times I)$ is {\it straight}, i. e. $f(S^{p-1}\times I)=f(S^{p-1}\times 0)\times I$.\label{proc-item3}
\end{enumerate}
\end{proc-definition}


\begin{procfigure}
\includegraphics{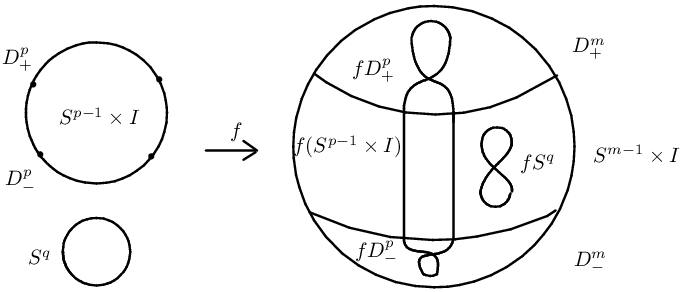}
\caption{A standardized link map}
\label{proc-fig1}
\end{procfigure}

\begin{proc-lemma}\label{proc-lem3}
Suppose that $p\le q+1$, $p\le m-3$ and $2p+2q\le
3m-5$\textup{;} then any link map $f:S^p\sqcup S^q\to S^m$ is link
homotopic to a standardized link map.
\end{proc-lemma}

\begin{proof}[Proof of the surjectivity in~\ref{proc-th2pr} for $p\le q$ modulo~\ref{proc-lem3}.] 
Take a link map $f:S^{p+1}\sqcup S^{q}\to S^{m+1}$. Let us modify it to a suspension
by a link homotopy.

By Lemma~\ref{proc-lem3} we may assume that $f$ is standardized.

Push $fS^q$ along the fibers of $S^{m}\times I$ until it lies in $S^{m}\times 0=\partial D^{m+1}_+$.
After that transform $fD^{p+1}_+$ and $f(S^{p+1}-\interior D^{p+1}_+)$ to the cones over
$f\partial D^{p+1}_+$ in $D^{m+1}_+$ and
in $S^{m+1}-\interior D^{m+1}_+$, respectively (by a rectilinear link homotopy).
The obtained link map is the
suspension of a link map $S^p\sqcup S^q\to S^m$.
\end{proof}

Now we proceed to the proof of Lemma~\ref{proc-lem3}. First we prove it for $p\le q$, then for
$p=q+1$. The proofs of all technical claims below can be skipped for the
first reading. From now till the end of \S\ref{proc-sect2} we work in piecewise linear category.

\begin{proof}[Proof of Lemma~\ref{proc-lem3} for $p\le q$] Let us make a generic link map $f:S^p\sqcup S^q\to S^m$
standardized by performing certain homeomorphisms~$S^p\to S^p$ and $S^m\to S^m$. (Formally, performing homeomorphisms
$h_p:S^p\to S^p$ and $h_m:S^m\to S^m$ means a link homotopy transforming $f$ to
$h_m\circ f\circ (h_p^{-1}\sqcup \identity_{S^q})$.)


\smallskip


\noindent(1) {\it Construction of the homeomorphism $S^p\to S^p$\textup{:} splitting the sphere $S^p$.}
(The Zeeman engulfing, see Figure~\ref{proc-fig2} to the left.)
Consider the self-intersection set of the $p$-component $S(f)=\closure\{\,x\in S^p:|f^{-1}fx|\ge2\,\}$.
Let $A_+$ be the skeleton of $S(f)$ formed by the simplices of
dimension not greater than $\frac{1}{2}\dim S(f)$ (in a triangulation of $S^p$, $S^q$ and $S^m$
such that $f:S^{p}\sqcup S^{q}\to S^{m}$ is simplicial). Let
$A_-$ be the subcomplex dual to $A_+$ (i. e., $A_-$ is the subcomplex formed by
all simplices $\sigma$ of the first barycentric subdivision of $S(f)$
such that $\sigma\cap A_+=\emptyset$).

\begin{proc-claim}\label{proc-cla}
There exist subpolyhedra $B_\pm\subset S^p$ such that $B_{\pm}\supset A_{\pm}$ and $B_\pm\cong CA_{\pm}$.
\end{proc-claim}

\begin{proof} Generically $\dim S(f)\le 2p-m$, so
$\dim A_{\pm}\le p-[\frac{m+1}{2}]$. Take generic
extensions $i_\pm:CA_\pm\to S^{p}$ of the inclusions $A_\pm\hookrightarrow
S^{p}$. They are embeddings, because $2(p-[\frac{m+1}{2}]+1)-p<0$
by the assumption
$p\le m-3$. Put $B_\pm=i_\pm CA_\pm$.
\end{proof}

\begin{proc-claim}\label{proc-clb} Generically $B_+\cap B_-=\emptyset$ and $B_{\pm}\cap
S(f)=A_{\pm}$.
\end{proc-claim}

\begin{proof}
This follows from \ 
$\dim(B_+\cap B_-)\le 2(p-[\frac{m+1}{2}]+1)-p<0$
and
\begin{equation*}
\dim(B_\pm-A_\pm)\cap S(f)\le(p-\left[\tfrac{m+1}{2}\right]+1)+(2p-m)-p\le
\tfrac{1}{2}(2q+2p-3m+4)<0,
\end{equation*}
which is a corollary of the assumptions $p\le q+1$ and $2p+2q\le 3m-5$.
\end{proof}

Take disjoint regular neighborhoods of $B_\pm$ in $S^{p}$. Perform
an orientation-preserving homeomorphism $S^p\to S^p$ taking these neighbourhoods to
the balls of the standard decomposition $S^p=D^p_+\cup (S^{p-1}\times I)\cup D^p_-$.

\begin{procfigure}
\includegraphics{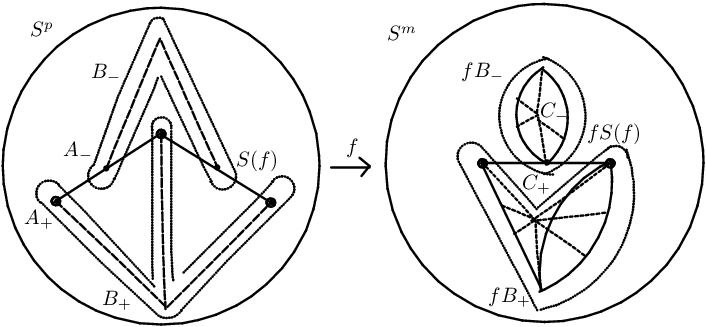}
\caption{Splitting of the spheres $S^p$ and $S^m$}
\label{proc-fig2}
\end{procfigure}



\smallskip

\noindent(2) {\it Construction of the 1st homeomorphism $S^m\to S^m$\textup{:} splitting the sphere $S^m$.}
(See Figure~\ref{proc-fig2} to the right.)
Take subpolyhedra $C_\pm\subset S^m$ such that $C_\pm\supset fB_\pm$ and $C_\pm\cong CfB_\pm$
(constructed analogously to Claim~\ref{proc-cla}).

\begin{proc-claim}\label{proc-clc}
Generically $C_+\cap C_-=\emptyset$, $C_\pm\cap fS^p=fB_\pm$ and $C_\pm\cap fS^q=\emptyset$.
\end{proc-claim}

\begin{proof}
This follows from the inequalities $2(p-[\frac{m+1}{2}]+2)-m<0$,
which holds by the assumption $p\le m-3$, and
$\dim(C_\pm-fB_\pm)\cap fS^{p}\le (p-[\tfrac{m+1}{2}]+2)+p-m<0$,
which holds by the assumptions $p\le q$ and $2p+2q\le 3m-5$, and
$\dim(C_+\cap fS^q)\le (p-[\frac{m+1}{2}]+2)+q -m <0$,
which is equivalent to $2p+2q\le 3m-5$.
\end{proof}

Take disjoint
regular neighborhoods of $C_\pm\cup fD^{p}_\pm$ relatively $f\partial D^{p}_\pm$ in
$S^{m}-fS^q$. Perform a homeomorphism $S^m\to S^m$
taking them to the balls of the standard
decomposition $S^m=D^m_+\cup (S^{m-1}\times I)\cup D^m_-$.
By Claim~\ref{proc-clc} the obtained link map satisfies properties~(\ref{proc-item1}) and~(\ref{proc-item2}) of a standardized link map (see Definition~\ref{proc-defslm}).

To satisfy property~(\ref{proc-item3}) perform the following homeomorphism $S^m\to S^m$.

\smallskip

\noindent(3) {\it Construction of the 2nd homeomorphism $S^{m}\to S^{m}$\textup{:} straightening $f(S^{p-1}\times I)$.}
\begin{proc-claim}\label{proc-cld}
There is a homeomorphism $S(f)\cap (S^{p-1}\times I)\cong (S(f)\cap (S^{p-1}\times 0))\times I$
taking $S(f)\cap (S^{p-1}\times j)$ to $(S(f)\cap (S^{p-1}\times 0))\times j$ for $j=0,1$.
\end{proc-claim}

\begin{proof} Take the first barycentric subdivision of the triangulation from step (1).
Then each simplex $\sigma\subset S(f)$ such that
$\sigma\not\subset A_+\cup A_-$ is the join of two simplices
$\sigma_+\subset A_+$ and $\sigma_-\subset A_-$. By Claim~\ref{proc-clb} the polyhedra $\sigma\cap
D^p_{\pm}$ are regular neighborhoods of $\sigma_{\pm}$ in
$\sigma$. So there is a natural homeomorphism $\sigma\cap
(S^{p-1}\times I)\cong(\sigma\cap(S^{p-1}\times 0))\times I$.
Combining such
homeomorphisms for all simplices $\sigma\subset S(f)$ such that
$\sigma\not\subset A_+\cup A_-$, we get the required homeomorphism $S(f)\cap (S^{p-1}\times
I)\cong (S(f)\cap (S^{p-1}\times 0))\times I$.
\end{proof}

\begin{proc-claim}\label{proc-cle} There is a homeomorphism $f(S^{p-1}\times I)\cong f(S^{p-1}\times 0)\times I$ taking
$f(S^{p-1}\times j)$ to $f(S^{p-1}\times 0)\times j$ for $j=0,1$.
\end{proc-claim}

\begin{proof} (The Alexander trick)
By Claim~\ref{proc-cld} the inclusion $i:S(f)\cap (S^{p-1}\times I)\hookrightarrow S^{p-1}\times I$
is a concordance. Perform an ambient isotopy
of $S^{p-1}\times I$ making $i$ an isotopy. Since any isotopy is ambient, there is
a homeomorphism $h:S^{p-1}\times I\to S^{p-1}\times I$ such that
$h(S(f)\cap (S^{p-1}\times I))$ is straight, i. e. is equal to $h(S(f)\cap (S^{p-1}\times 0))\times I$.
The required homeomorphism $f(S^{p-1}\times I)\cong f(S^{p-1}\times 0)\times I$ is the quotient of $h$.
(Analogously to the proof of~\ref{proc-cld} it can be checked that this quotient is well-defined.)
\end{proof}

\begin{proc-claim}\label{proc-clf} There is a homeomorphism of $S^{m-1}\times I$
making $f(S^{p-1}\times I)$ straight.
\end{proc-claim}

\begin{proof}By Claim~\ref{proc-cle} the inclusion $i:f(S^{p-1}\times I)\hookrightarrow S^{m-1}\times I$
is a concordance. Arguing as in the proof of Claim~\ref{proc-cle} we get the required homeomorphism.
\end{proof}

Perform a homeomorphism $S^m\to S^m$ extending the one given by Claim~\ref{proc-clf}. The obtained link map
is standardized, so for $p\le q$ Lemma~\ref{proc-lem3} is proved.
\end{proof}


%


\begin{proof}[Proof of Lemma~\ref{proc-lem3} for $p=q+1$] The proof
is analogous to the proof in case $p\le q$, only the cones $B_\pm$ and $C_\pm$ should be replaced
by collapsible polyhedra given by the following claim, cf. \cite{Hud69}.
\end{proof}

\begin{proc-claim}\label{proc-clg}
There are collapsible subpolyhedra $B_{\pm}\subset S^p$ and
$C_{\pm}\subset S^m$ satisfying Claims~\ref{proc-clb} and~\ref{proc-clc}.
\end{proc-claim}


\begin{proof} (The Irwin trick) Let $\bar B_\pm$ be the polyhedra
given by Claim~\ref{proc-cla}. Define $\bar C_\pm$ analogously. These polyhedra satisfy all
the required properties except $\bar C_\pm\cap fS^p=f\bar B_\pm$.
By the inequality
$\dim(\bar C_\pm-f\bar B_\pm)\cap fS^{p}\le (p-[\tfrac{m+1}{2}]+2)+p-m\le 0$
the set $(\bar C_\pm-f\bar B_\pm)\cap fS^p$ (if nonempty) consists of finitely many points
not belonging to $S(f)$. Join each of these points with $\bar B_\pm$ by
a generic arc in $S^{p}$. Let $B_\pm$ be the union of these arcs
and the cone $\bar B_\pm$. Adding appropriate cones over $f(B_\pm-\bar B_\pm)$ to
$\bar C_\pm$, we get a collapsible polyhedron $C_\pm\subset S^{m}$ such that $\dim(C_\pm-\bar C_\pm)\le 2$.
The polyhedra $B_{\pm}$ and $C_{\pm}$ are the required.
\end{proof}



The injectivity in Theorem~\ref{proc-th2pr} is proved by a relative version of the above argument.
A {\it standardized} link map $f:D^p\sqcup D^q\to D^m$
is defined as in~\ref{proc-defslm}, only we fix the standard decomposition of the {\it disc}
$D^k=D^k_+\cup (D^{k-1}\times I)\cup D^k_-$
instead of the {\it sphere}.
Denote by $D^{k-1}_\pm=D^{k}_\pm\cap\partial D^k$. Assume that $\partial D^{k-1}_+$ is the equator of $\partial D^k$.

\begin{proc-lemma}\label{proc-lem3pr}
Suppose that $p\le q+1$, $p\le m-3$ and $2p+2q\le
3m-5$\textup{.} Then any generic proper link map $f:D^p\sqcup D^q\to D^m$, whose
restriction to the boundary is a suspension, is link homotopic
(relatively the boundary) to a standardized link map.
\end{proc-lemma}

\begin{proof}[Proof of the injectivity in~\ref{proc-th2pr} for $p\le q$ modulo~\ref{proc-lem3pr}]
By 'link concordance implies link homotopy' theorem \cite{Mel00}, it suffices to prove that if the suspension of a link map
$f_0:S^p\sqcup S^q\to S^m$ is null link concordant, then the link map $f_0$
is null link concordant. Take a null link concordance $f:D^{p+2}\sqcup
D^{q+1}\to D^{m+2}$ of $\Sigma f_0$. By Lemma~\ref{proc-lem3pr} we may assume that
the link map $f$ is standardized.

Push $fD^{q+1}$ along the fibers of $D^{m+1}\times I$
toward $\partial D^{m+2}_+$ until it lies in $\partial D^{m+2}_+ -
\partial D^{m+2}$. The restriction of the obtained link map to
$\closure(\partial D^{p+2}_+ - \partial D^{p+2})\sqcup D^{q+1}$ is the required null
link concordance of the link map $f_0$.
\end{proof}


\begin{proof}[Proof of Lemma~\ref{proc-lem3pr}] The proof is analogous to the proof of Lemma~\ref{proc-lem3} with the following
modifications. 
Let $\hat D^k$ be the ball obtained from $D^k$ by attaching two cones $C D^{k-1}_\pm$ along $D^{k-1}_\pm$. Let
$\hat f:\hat D^p\sqcup \hat D^q\to \hat D^m$ be the obvious extension
of the link map $f:D^p\sqcup D^q\to D^m$.
Clearly, it suffices to make the link map $\hat f$ standardized (performing homeomorphisms of $\hat D^p$ and $\hat D^m$ fixed on the boundary).

\smallskip
\noindent
(1) {\it  Construction of a homeomorphism $\hat D^p\to\hat D^p$ for $p\le q$.}
Let $A_+$ be the union of $\partial D^{p-1}_+$ and
all simplices of $S(f)$ having the dimension not greater than $\frac{1}{2}\dim{S(f)}$.
Let $A_-$ be the subcomplex dual to $A_+$.

\begin{proc-claim}\label{proc-clh}
There are subpolyhedra $B_\pm\subset D^p$
collapsible to $B_\pm\cap D^{p-1}_\pm$ and satisfying Claim~\ref{proc-clb}.
\end{proc-claim}

\begin{proof}
Take a generic homotopy $i_t:A_\pm\to D^p$ fixed on
$A_\pm\cap D^{p-1}_\pm$, such that $i_0:A_\pm\hookrightarrow D^p$ is the
inclusion and $i_1A_\pm\subset D^{p-1}_\pm$. Let $B_\pm$ be the
trace of~$i_t$.
\end{proof}

Take appropriate regular neighborhoods of $B_\pm\cup C(B_\pm\cap D^{p-1}_\pm)$
in $D^p\cup CD^{p-1}_\pm$.
Perform a homeomorphism $\hat D^p\to\hat D^p$ fixed on the boundary, taking them to the balls
of the standard decomposition of $\hat D^p$.

\smallskip
Steps~(2) and~(3) from the proof of Lemma~\ref{proc-lem3} are modified analogously.
\end{proof}

Thus we have proved Theorem~\ref{proc-th2pr} for $p\le q$.

\begin{proof}[\it Proof of Theorem~\ref{proc-th2pr} for $p>q$] The map $\Sigma$ is surjective
as the composition
\begin{equation*}
\xymatrix@1{
{{LM}^m_{p,q}}             \ar[rr]^{\Sigma^{p-q}}        &&
{{LM}^{m+p-q}_{p,p}}       \ar[rr]^{\Sigma}              &&
{{LM}^{m+p-q+1}_{p+1,p}}   \ar[rr]^-{(\Sigma^{p-q})^{-1}} &&
{{LM}^{m+1}_{p+1,q}},
}
\end{equation*}
in which all the maps are well-defined and surjective by the case $p\le q$ of Theorem~\ref{proc-th2pr}.
The injectivity of $\Sigma$ is proved analogously.
\end{proof}

\procsection{Classification of links}\label{proc-sect3}

We prove Theorem~\ref{proc-th1} as follows. First we prove a suspension theorem
for links (Lemma~\ref{proc-th4}) reducing
the classification of links to the classification of {\it disc link maps}.
Then we simplify the group of disc link maps, and find it using the
classification of link maps. Formally,
\ref{proc-th1} follows from~\ref{proc-th2pr}, \ref{proc-cor4}, \ref{proc-lem5}, \ref{proc-th6} and 5-lemma.

Let us introduce some notation. Throughout \S\ref{proc-sect3} we work in smooth category.

Denote by $\widehat{L}^m_{p,q}$ the
group of concordance classes of embeddings $S^p\sqcup S^q\to S^m$, whose restrictions both to $S^p$ and to $S^q$ are unknotted.

An {\it almost link} is a link map $f:S^p\sqcup S^q\to S^m$ whose
restriction to $S^q$ is an unknotted embedding. An {\it almost concordance} is a
link concordance $(S^p\sqcup S^q)\times I\to S^m\times I$, whose restriction to $S^q\times I$
is a concordance.
Let $\overline{L}^m_{p,q}$ be the set of almost links
up to almost concordance.
For $p,q\le m-3$ this set is a commutative group with respect to 'componentwise connected sum' operation.
It is not difficult to see that this group is isomorphic to $\pi_p(S^{m-q-1})$ (cf. Definition of $\lambda$ below).


A {\it disc link map} is a proper link map
$f:D^p\sqcup D^q\to D^m$ whose restriction to $(\partial D^p)\sqcup D^q$ is an embedding and $f:\partial D^p\to \partial D^m$ is unknotted. A {\it disc link concordance} is a proper link
concordance $(D^p \sqcup D^q)\times I\to D^m\times I$, whose restriction to $(\partial D^p\sqcup D^q)\times I$ is a
concordance. Let $\widehat{DM}^m_{p,q}$ be the set of disc link maps up
to disc link concordance. For $p,q\le m-3$ it has a natural commutative group
structure.

The following two results reduce the classification of links to the classification of disc link maps.

\begin{proc-lemma}[Geometric EHP sequence for links]\label{proc-th4}
\textup{(A. Skopenkov, cf. \cite{Ne84}, \cite[Theorem~1.9]{Hae66A}, see Figure~\ref{proc-fig3})} For $p,q\le m-3$
there is an exact sequence\textup{:}
\begin{equation*}
\xymatrix@1{
{\dots}                       \ar[r]       &
{\widehat{L}^m_{p,q}}         \ar[r]^{e}   &
{\overline{L}^m_{p,q}}        \ar[r]^{h}   &
{\widehat{DM}^m_{p,q}}        \ar[r]^{p}   &
{\widehat{L}^{m-1}_{p-1,q-1}} \ar[r]       &
{\dots}
}
\end{equation*}
\end{proc-lemma}


\begin{procfigure}
\includegraphics{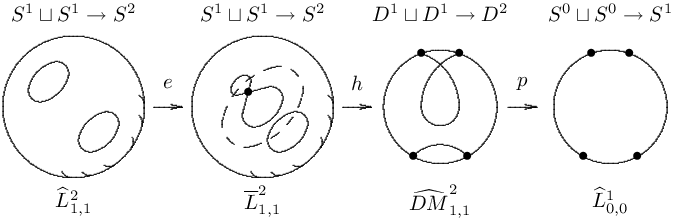}
\caption{Geometric EHP sequence for links}
\label{proc-fig3}
\end{procfigure}

\begin{proof}
{\it Construction of the homomorphisms.}
Let $e$ be the obvious map.
Let $p$ be the 'restriction to the boundary' map. It is well-defined because $\partial D^q\to \partial D^m$ extends to an embedding $D^q\to D^m$ and hence is unknotted.
Let $h$ be the 'cutting' homomorphism defined as
follows. Take a generic almost link $f:S^p\sqcup S^q\to S^m$.
Take generic points $x\in fS^p$ and
$y\in fS^q$. Join them by a path $l$ meeting
$f(S^p\sqcup S^q)$ only at $\partial l$. Let $\bar D^m$ be the
complement to a small neighborhood of $l$ in $S^m$. Denote by
$\bar D^p\sqcup\bar D^q=f^{-1}\bar D^m$. Set $h(f)$ to be the
restriction of $f$ to a map $\bar D^p\sqcup\bar D^q\to \bar D^m$.

{\it Proof of the exactness.} We have
$\image p=\kernel e$ because a link $f:S^p\sqcup
S^q\to S^m$ extends to a disc link map $D^{p+1}\sqcup D^{q+1}\to
D^{m+1}$ if and only if it is almost concordant to a trivial link.
We have $\image h=\kernel p$ because a disc link
map $f:D^p\sqcup D^q\to D^m$ extends without adding new
self-intersections to an almost link $S^p\sqcup S^q\to S^m$ if and only if the restriction of $f$ to the boundary is null-concordant.


To prove $\image e\subset \kernel h$, take a proper embedding $f:D^p\sqcup D^q\to D^m$.
Take a pair of points
$x\in D^p$ and $y\in D^q$. Join $fx$ and $fy$ by an arc $l$
meeting $f(D^p\sqcup D^q)$ only at $\partial l$. Let $\bar D^m$ be
a small neighborhood of $l$ in $S^m$. Denote by
$\bar D^p\sqcup\bar D^q=f^{-1}\bar D^m$.
The restriction $f:(D^p-\bar D^p)\sqcup (D^q-\bar D^q)\to (D^m-\bar D^m)$ is a concordance.
By 'concordance implies isotopy' theorem
we may assume that this restriction is level-preserving.
Then the Alexander trick shows that the embedding $f:D^p\sqcup D^q\to D^m$ is ambient isotopic to the
restriction $f:\bar D^p\sqcup \bar D^q\to \bar D^m$. The latter embedding is trivial, thus $h\circ e=0$.

To prove $\kernel h\subset\image e$, take $f\in
\overline{L}^m_{p,q}$ such that $h(f)=0$. By definition, there
exist a disc link concordance $c$ between $h(f)$ and an
embedding. 
We may assume that the restriction of $c$ to the boundary is an isotopy. By isotopy extension theorem \cite{Hud69} it
extends to an ambient isotopy of the disc $S^m-\bar D^m$ (from the above definition of the map $h$). So $c$ can be extended to an almost concordance between $f$ and a link $f'\in \widehat{L}^m_{p,q}$. Hence
$f=e(f')$.
\end{proof}

\begin{proc-corollary}\label{proc-cor4}
$L^m_{p,q}\cong C^{m-p}_p\oplus C^{m-q}_q\oplus \pi_p(S^{m-q-1})\oplus
\widehat{DM}^{m+1}_{p+1,q+1}$ for $p\le q\le m-3$.
\end{proc-corollary}

\begin{proof}
By \cite[Th. 2.4]{Hae66C} we have $L^m_{p,q}\cong
\widehat{L}^m_{p,q}\oplus C^{m-p}_p\oplus C^{m-q}_q$. So it suffices to show
that for $p\le q$ the homomorphism $e:\widehat{L}^m_{p,q}\to
\overline{L}^m_{p,q}$ 
has a right inverse. 
The required right
inverse $\pi_p(S^{m-q-1})\to \widehat{L}^m_{p,q}$ was constructed in \cite[Th. 10.1]{Hae66C}: it
takes the class
of a map $\phi:S^p\to S^{m-q-1}$ to a link $f:S^p\sqcup S^q\to
D^{q+1}\times S^{m-q-1}\subset S^m$ given by
the formula $f(x\sqcup y)=(\frac{1}{2}x;\phi x)\sqcup (y;c)$, where $c\in S^{m-q-1}$ is fixed.
\end{proof}

Let us simplify the group $\widehat{DM}^{m}_{p,q}$. Define $\overline{DM}^m_{p,q}$ to be the group of proper
link maps $f:D^p\sqcup D^q\to D^m$ whose restriction $f:\partial
D^p\to \partial D^m$ is an unknotted embedding (up to link concordance whose restriction to
$\partial D^p\times I$ is a concordance).

\begin{proc-lemma}\label{proc-lem5}
If $p,q\le m-3$ then the natural map
$\widehat{DM}^m_{p,q}\to\overline{DM}^m_{p,q}$ is bijective for
$2p+2q\le 3m-5$ and surjective for $2p+2q\le 3m-4$.
\end{proc-lemma}

\begin{proof}{\it The surjectivity.} Take a generic link map
$f\in \overline{DM}^m_{p,q}$. The pair $(D^m-fD^p,\partial D^m-f\partial D^p)$ is $(2m-2p-3)$-connected,
because $H_i(D^m-fD^p,\partial D^m-f\partial D^p) \cong
H^{m-i-1}(fD^p)=0$ for $i\le 2m-2p-3$ (because $fD^p$ is homotopy equivalent to the mapping cone of the restriction $f:S(f)\to fS(f)$,
having the dimension at most $2p-m+1$, cf. \cite[Lemma~4.2]{HaKa98}). Thus by the assumptions $q\le m-3$, $2p+2q\le 3m-4$ and
the embedding theorem moving the boundary \cite{Hud69} the restriction
$f\left|_{D^q}\right.:(D^q,\partial D^q)\to (D^m-fD^p,\partial
D^m-f\partial D^p)$ is homotopic to an embedding. So $f$ belongs to the image of the natural map
$\widehat{DM}^m_{p,q}\to\overline{DM}^m_{p,q}$.

{\it The injectivity.} Take a generic link concordance $f:(D^p\sqcup D^q)\times I\to D^m\times I$, whose
restriction to $D^q\times\partial I\cup
\partial D^p\times I$ is an embedding.
It suffices to remove the
self-intersection of $D^q\times I$ by a link homotopy fixed on
$(D^p\sqcup D^q)\times \partial I$. It is possible by
the following theorem proved similarly to \cite{Hud69}, because the
pair $(D^m\times I-f(D^p\times I),\partial D^m\times I-f(\partial
D^p\times I))$ is $(2m-2p-3)$-connected.
\end{proof}

\begin{proc-theorem}[Embedding theorem moving a part of the boundary] \label{proc-embth}
Let $M^{m+1}$\textup{,} $Y^{m}\subset\partial M$ and $X^q\subset\partial D^{q+1}$ be compact manifolds. Let $f:(D^{q+1},X)\to (M,Y)$ be a proper map such that $f\left|_{\partial D^{q+1}-X}\right.$ is an embedding. If $q\le
m-3$ and $(M;Y)$ is $(2q-m+2)$-connected, then $f$ is properly
homotopic $\relatively\partial D^{q+1}-X$ to an embedding.
\end{proc-theorem}

Let us find the group $\overline{DM}^m_{p,q}$. Denote by $n=p+q+1-m$.
We are going to define a homomorphism $\beta:\overline{DM}^m_{p,q}\to \pi_{n}(V_{M+m-p-1,M})$.
The following theorem and 5-lemma imply the bijectivity of this homomorphism.

\begin{proc-theorem}\label{proc-th6}
\textup{(A. Skopenkov, cf. \cite[Th. 3.1]{Kos90}, \cite[Lemma 5.1]{Ker59}, \cite[Th. 4.8]{Kos88})} For $p,q\le m-3$ and $3p+q\le 3m-5$ there is the following
diagram with exact lines, commutative up to sign\textup{:}
\begin{equation*}
\xymatrix{
{\overline{L}^m_{q,p}}         \ar[r]^{e} \ar[d]^{\lambda}  & {LM^m_{p,q}}                    \ar[r]^{h}  \ar[d]^{\alpha}  &
{\overline{DM}^m_{p,q}}        \ar[r]^{p} \ar[d]^{\beta}    & {\overline{L}^{m-1}_{q-1,p-1}}  \ar[r]      \ar[d]^{\lambda} & {\dots} \\
{\pi_q(S^{m-p-1})}             \ar[r]^-{E}                  & {\pi^S_{n}}                     \ar[r]^-{H}                  &
{\pi_{n}(V_{M+m-p-1,M})}       \ar[r]^-{P}                  & {\pi_{q-1}(S^{m-p-1})}          \ar[r]                       & {\dots}
}
\end{equation*}
\end{proc-theorem}






\begin{proof} The top line is defined analogously to Lemma~\ref{proc-th4} (with
similar proof of the exactness). The bottom line is the stable James
EHP sequence, for which we use the following geometric construction \cite[\S1 and \S4]{KoSa77}, cf. \cite{Jam54, Szu76, Ecc80}.

{\it Construction of the EHP sequence.}
Identify the groups $\pi_q(S^{m-p-1})$ and $\pi^S_{n}$ with the groups of framed embeddings and
immersions, respectively, of closed $n$-manifolds into $S^q$ (up to framed cobordism).
A {\it proper immersion} is a proper framed immersion of an
$n$-manifold into $D^q$, whose restriction to the boundary is an
embedding. A {\it proper cobordism} is a proper framed immersion
$c:N^{n+1}\to D^q\times I$, whose restriction to
$c^{-1}(S^{q-1}\times I)$ is an embedding. By \cite[Prop.~4.1]{KoSa77} we can identify $\pi_{n}(V_{M+m-p-1,M})$ with
the group of proper immersions up to proper cobordism.

Let $E:\pi_q(S^{m-p-1})\to \pi^S_{n}$ be the obvious map and let
$P:\pi_{n}(V_{M+m-p-1,M})\to\pi_q(S^{m-p-1})$ be 'restriction to the boundary' map.
Let $H:\pi^S_{n}\to\pi_{n}(V_{M+m-p-1,M})$ be cutting homomorphism, defined by removing
small discs from an immersed $n$-manifold and the sphere $S^q$.

%

{\it Construction of the vertical homomorphisms.} Remove a point
from $S^m$ and identify the result with $\Bbb R^m$. For a link map
$f:X\sqcup Y\to \Bbb R^m$ define the map $\tilde f:X\times Y\to
S^{m-1}$ by the formula
$
\tilde f(x,y)=\frac{fx-fy}{|fx-fy|}.
$
Denote by $\projection:X\times Y\to Y$ the
obvious projection.

{\it Definition of $\alpha$.} (cf. \cite{Kos88}) Let $f:S^p\sqcup S^q\to \Bbb R^m$ be a
generic link map. By general position $v=(1,0,\dots,0)$ is a regular value of~$\tilde f$. The Pontryagin--Thom construction gives a framed embedding $\tilde f^{-1}v\to S^p\times S^q$. Perform a regular homotopy to get a framed immersion $c\colon \tilde f^{-1}v\to S^p\times S^q$ with the last $p$ vectors of the framing tangent to the $S^p$-fibers of the product $S^p\times S^q$ and forming a positive basis of the tangent space. (Such immersion exists by Hirsch theory; see \cite[Theorem~6.1 and~5.9]{Hir59}). Let $\alpha(f)$ be the cobordism class of the framed immersion
$\projection\circ c\colon \tilde f^{-1}v\to S^q$. The framing is given by the projection of the first $q-n$ vectors of the framing of $c$.
(Clearly, $\alpha$ commutes with $\Sigma$, hence by~\ref{proc-th2pr} and \cite[Th.~ 2.13]{Kos88} $\alpha$ is an isomorphism
for $p,q\le m-3$ and $2p+2q\le 3m-5$.)

{\it Definition of $\lambda$.} (cf. \cite{Hae66C}) Take a generic link map $f\in \bar L^m_{q,p}$ whose restriction to $S^p$ is an unknotted embedding. The complement $S^m-fS^p$
retracts to a sphere $S^{m-p-1}$ bounding a normal disc to $fS^p$.
Put the image $fS^q$ into the sphere $S^{m-p-1}$ by an appropriate link homotopy fixed on $S^p$. By an isotopy put $fS^p\sqcup S^{m-p-1}$ onto the \emph{standard link} consisting of two spheres
\begin{align*}
& (x_{p+1}-1)^2+ x_{p+2}^2+\dots +x_{m}^2=1, \quad x_{1}=\dots=x_{p}=0 \text{ and } \\
& x_{1}^2+\dots+x_{p+1}^2=1, \quad x_{p+2}=\dots=x_{m}=0.
\end{align*}
Fix an orientation of the former sphere and assume that $f\left|_{S^p}\right.$ preserves orientation. By general position $v=(1,0,\dots,0)$ is a regular value of the map $\tilde f$. Let $\lambda(f)$ be the cobordism class of the framed embedding $\projection\colon\tilde f^{-1}v\to S^q$.

{\it Proof that $\lambda$ is well-defined}. Let us prove that $\projection\colon\tilde f^{-1}v\to S^q$ is indeed a framed embedding. Moreover, we show that $\tilde f^{-1}v=f^{-1}(0,\dots,0)\times f^{-1}(1,0,\dots,0)$ is a framed submanifold of $f^{-1}(0,\dots,0)\times S^q$.

Take any $x\in S^p$ and $y\in S^q$ such that $\tilde f(x,y)=v$. Since $fy\in fS^q\subset S^{m-p-1}$ and $v$ is parallel to the subspace $\Bbb R^{m-p}\supset S^{m-p-1}$, it follows that $fx\in\Bbb R^{m-p}$. The subspace intersects $fS^p$ at the origin and the point with all the coordinates vanishing except $x_{p+1}=2$. Since the ray parallel to $v$ starting at the latter point does not intersect $S^{m-p-1}$, it follows that $x=f^{-1}(0,0,\dots,0)$. Since the ray parallel to $v$ starting at the origin intersects $S^{m-p-1}$ only at $(1,0,\dots,0)$, it follows that $y\in f^{-1}(1,0,\dots,0)$. Thus $\tilde f^{-1}v=f^{-1}(0,\dots,0)\times f^{-1}(1,0,\dots,0)$ and $\projection\colon\tilde f^{-1}v\to S^q$ is an embedding. Since the last $p$ basis vectors of $\mathbb{R}^m$ are tangent to $fS^p$ at the origin and orthogonal to $S^{m-p-1}$ at $(1,0,\dots,0)$, it follows that the last $p$ vectors of the framing of $\tilde f^{-1}v$ are tangent to the $S^p$-fiber. Thus the framing of $\projection\tilde f^{-1}v$ is well-defined and coincides with the framing of $f^{-1}(1,0,\dots,0)$.

Independence $\lambda (f)$ on the choice of the isotopy taking $fS^p\sqcup S^{m-p-1}$ onto the standard link is checked analogously to \cite[Th.~7.1]{Hae66C}.

{\it Definition of $\beta$.} (cf. \cite{Kos90}) The definition of $\beta$ is a relative version of the above two ones: at the boundary the construction is analogous to $\lambda$, in the interior --- to $\alpha$. Take a generic proper disc link map
$f:D^p\sqcup D^q\to \Bbb R^m_+$, where $\Bbb R^m_+$ is the upper
half-space $x_m\ge 0$. By a proper link homotopy, restricting to an isotopy of
$\partial D^p$, put the images $f\partial D^p$ and $f\partial D^q$ into the components of the
standard link in $\partial\Bbb R^m_+$. By general position $v=(1,0,\dots,0)$ is a regular value of $\tilde f$.
We get a proper framed embedding $\tilde f^{-1}v\to D^p\times D^q$. The last $p$ vectors of the framing at the boundary are tangent to the $D^p$-fibers (as shown in the next paragraph). Perform a regular homotopy fixed on the boundary to get a proper immersion $c\colon \tilde f^{-1}v\to D^p\times D^q$ with the last $p$ vectors of the framing of the interior tangent to the $D^p$-fibers as well. Let $\beta(f)$ be the proper cobordism class of
$\projection\circ c\colon \tilde f^{-1}v\to D^q$.

{\it Proof that $\beta$ is well-defined}. Let us prove that $\projection\circ c\colon \tilde f^{-1}v\to D^q$ is indeed a proper immersion.
Moreover, we show that the restriction of the map to the boundary coincides with $\lambda(f\left|_{\partial D^p\sqcup\partial D^q}\right.)$ up to sign.

Since $v$ is parallel to $\partial \Bbb R^m_+$ and $f(\interior D^p),f(\interior D^q)\subset \interior \Bbb R^m_+$, it follows that $\tilde{f}^{-1}v$  intersects neither $\interior D^p \times\partial D^q$ nor $\partial D^p \times \interior D^q$. Thus $\partial\tilde{f}^{-1}v\subset\partial D^p \times\partial D^q$ and $\projection\colon \partial\tilde f^{-1}v\to \partial D^q$ coincides with $\lambda(f\left|_{\partial D^p\sqcup\partial D^q}\right.)$ up to sign. As shown in the 'proof that $\lambda$ is well-defined' above, $\lambda(f\left|_{\partial D^p\sqcup\partial D^q}\right.)$ is an embedding and the last $p-1$ vectors of the framing of the embedding $\partial\tilde f^{-1}v\to \partial D^q\times \partial D^q$  are tangent to the $D^p$-fibers. 

Thus the required diagram is constructed. The commutativity up to sign is checked directly. Let us do it just for one square.

{\it The commutativity of the left square.} Take a generic link map $f\in \bar L^m_{q,p}$. Put the images $fS^p$ and $fS^q$ into the standard link by an almost concordance. Then $\lambda(f)$ and $E\lambda(f)$ are the classes of the framed embedding $\projection\colon\tilde f^{-1}v\to S^q$ up to embedded and immersed cobordism respectively. Since the last $p$ vectors of the framing of $\tilde f^{-1}v$ are already tangent to the $S^p$-fiber, it follows that $\alpha e(f)$ is the immersed cobordism class of $\projection\colon \tilde f^{-1}v\to S^q$. Thus $E\lambda=\pm\alpha e$.
\end{proof}


%
%



It is interesting to learn, if the diagrams of Theorems~\ref{proc-th6} and~\cite[Th. 3.1]{Kos90} are isomorphic.


\begin{proc-remark}\label{proc-rem1} The restriction $2p+2q\le 3m-6$ in Theorem~\ref{proc-th1} is best
possible, the formula fails for $2p+2q=3m-5$. For example, take
$m=p+4=4k-1$, $k\ge 5$, $q=2k+1$. Then the group
$LM^{m+1}_{p+1,q+1}$ is infinite \cite[p. 755-756]{Kos90}. Thus by \ref{proc-lem5} and ~\ref{proc-th6}
$\rank\widehat{DM}^{m+1}_{p+1,q+1}\ge\rank\overline{DM}^{m+1}_{p+1,q+1}>\rank\pi_{p+q+2-m}(V_{M+m-p-1,M})$.
So by \ref{proc-cor4} the rank of the left-hand side in the formula of~\ref{proc-th1} is greater than the rank
of the right-hand side.
\end{proc-remark}

\begin{proc-remark}\label{proc-rem2} The argument of \cite{Hae66C} can be extended to prove Theorem~\ref{proc-th1} at least for $2p_1+2p_2\le 3m-7$ (here we use the notation of \cite{Hae66C}).  Indeed, the only step of that proof, in which this restriction is not sufficient, is
\cite[Proposition~10.2]{Hae66C}. Since the group $\Lambda^{(q)}_{p_1}$ (respectively, $\Pi^{(q)}_{m-2}$) is generated by $\theta_k(i_1,i_2)$ (respectively, $[[i_2,i_1],i_2]$ and $\theta_{k+1}(i_1,i_2)$) for all $k\ge 0$ such that $kp_1+p_2\ge (k+1)(m-2)$, the proposition
follows. A possible reason why this improvement was not noticed in \cite{Hae66C} was
that the restriction $2p_1+2p_2\le 3m-7$ did not appear there in contrast to $3p_1+p_2\le 3m-7$.
\end{proc-remark}



\procsubsection{Acknowledgements} The author is grateful to A.~Skopenkov for constant attention to this work and to P.~Akhmetiev, S.~Avvakumov, S.~Melikhov and the anonimous referee for useful suggestions.


\bibliographystyle{amsplain}


\begin{procthebibliography}{99}

\bibitem{BaTe99} A. Bartels, P. Teichner, \textit{All two dimensional links are
null homotopic}, Geom. Topol. \textbf{3} (1999), p. 235--252.


\bibitem{CRS07} M. Cencelj, D. Repovs, M. Skopenkov, \textit{Homotopy type of
the complement to an immersion and classification of embeddings of tori},
Rus. Math. Surv. \textbf{62:5} (2007), p. 985--987,
\href{http://arxiv.org/abs/0803.4285v1}{arXiv:0803.4285v1}[math.GT].

\bibitem{CRS} M. Cencelj, D. Repov\v s, M. Skopenkov, \textit{Classification of knotted tori in 2-metastable dimension}, Mat. Sbornik \textbf{203:11} (2012), 129-158 (in Russian). English transl.: Sbornik Math. \textbf{203:11} (2012), p. 1654--1681; \href{http://arxiv.org/abs/0811.2745}{arXiv:math/0811.2745}[math.GT].

\bibitem{CFS14} D. Crowley, S. Ferry, M. Skopenkov, \textit{The rational classification of links in codimension $>2$}, Forum Math. \textbf{26:1} (2014), p. 239--269; \href{http://arxiv.org/abs/1106.1455}{arXiv:math/0811.2745}[math.GT].


\bibitem{Ecc80} P. Eccles, \textit{Multiple points of codimension one immersions},
Lect. Notes Math. \textbf{788} (1980), p. 23--38.

\bibitem{Ha86} N. Habegger, \textit{Knots and links in codimension greater than 2},
Topol. \textbf{25:3} (1986), p. 253--260.

\bibitem{HaKa98} N. Habegger, U. Kaiser, \textit{Link homotopy in the
2-metastable range}, Topol. \textbf{37:1} (1998), p. 75--94.

\bibitem{Hae66A} A. Haefliger, \textit{Differentiable embeddings of $S^n$ in
$S^{n+q}$ for $q>2$}, Ann. Math., Ser.3 \textbf{83} (1966) p. 402--436.

\bibitem{Hae66C} \bysame, \textit{Enlacements de spheres en codimension
superiure a 2}, Comm. Math. Helv. \textbf{41} (1966-67), p. 51--72 (in
French).

\bibitem{Hir59} M. W. Hirsch, \textit{Immersions of manifolds}, Trans. Amer. Math. Soc. \textbf{93:2} (1959), p. 242--276.

\bibitem{Hud69} J. F. P. Hudson, \textit{Piecewise-linear topology}, Benjamin, New
York-Amsterdam 1969.


\bibitem{Jam54} I. M. James, \textit{On the iterated suspension}, Quart. J. Math.
Oxford \textbf{5} (1954), p. 1--10.

\bibitem{Ker59} M. Kervaire, \textit{An interpretation of G. Whitehead's
generalization of H. Hopf's invariant}, Ann. Math. \textbf{69} (1959),
p. 345--362.

\bibitem{Kos88} U. Koschorke, \textit{Link maps and the geometry of their
invariants}, Manuscripta Math. \textbf{61:4} (1988), p. 383--415.


\bibitem{Kos90} \bysame, \textit{On link maps and their homotopy
classification}, Math.Ann. \textbf{286:4} (1990), p.753--782.

\bibitem{Kos97} \bysame, \textit{A generalization of Milnor's
$\mu$-invariants to higher dimensional link maps}, Topology \textbf{36:2}
(1997), p. 301--324.

\bibitem{KoSa77} U. Koschorke, B. Sanderson, \textit{Geometric interpretation
of the generalized Hopf invariant}, Math. Scand. \textbf{41} (1977), p.
199--217.

\bibitem{KrTa97} V. Krushkal, P. Teichner, \textit{Alexander duality, gropes
and link homotopy,} Geom. Topol. \textbf{1} (1997), p. 51--69.

\bibitem{Mel00} S. Melikhov, \textit{Link concordance implies link homotopy in
codimension $\ge3$}, Uspekhi Mat. Nauk \textbf{55:3} (2000), p. 183--184 (in
Russian). 

\bibitem{Mel} \bysame, \textit{Link concordance implies link homotopy in
codimension $\ge3$}, preprint.

\bibitem{Ne84} V. Nezhinsky, \textit{A suspension sequence in link theory}, Izv.
Akad. Nauk \textbf{48:1} (1984), p. 126--143. 

\bibitem{Pae56} G.F. Paechter, \textit{The groups $\pi_r(V_{n,m})$}, Quart. J.
Math. Oxford, Ser. 2, \textbf{7} (1956), p. 249--268.

\bibitem{RS99} D. Repovs and A. Skopenkov, \textit{New results on embeddings
of polyhedra and manifolds into Euclidean spaces}, Uspekhi Mat. Nauk \textbf{54:6} (1999), p. 61--109 (in Russian).

\bibitem{Sco68} G. P. Scott, \textit{Homotopy links}, Abh. Math. Sem. Univ. Hamburg
\textbf{32} (1968), p. 186--190.

\bibitem{S} A. Skopenkov, \textit{High codimension embeddings: classification}, submitted to Bull. Man. Atl.\\
{\tiny http://www.map.mpim-bonn.mpg.de/Embeddings\_in\_Euclidean\_space:\_an\_introduction\_to\_their\_classification}

\bibitem{M} A. Skopenkov, \textit{High codimension links}, submitted to Bull. Man. Atl.\\ \url{http://www.map.mpim-bonn.mpg.de/High_codimension_links},

\bibitem{Sko06} A. Skopenkov, \textit{Classification of knotted tori},
submitted, \href{http://arxiv.org/abs/1502.04470}{arXiv:math/1502.04470}[math.GT].

\bibitem{Sko07} \bysame, \textit{Embedding and knotting
of manifolds in Euclidean spaces}, \textit{in}: \textit{Surveys in Contemporary Mathematics}, \textit{Ed. N. Young and Y. Choi},
London Math. Soc. Lect. Notes \textbf{347} (2007), p. 248--342, \href{http://arxiv.org/abs/math.GT/0604045}{arXiv:math/0604045}[math.GT].

\bibitem{Szu76} A. Sz\"ucs, \textit{Cobordism group of $l$-immersions}, Acta Math.
Hungar. \textbf{28} (1976), p. 93-102. 

\end{procthebibliography}

\end{document}